# A stochastic power management strategy with skid avoidance for improving energy efficiency of in-wheel motor electric vehicles

**Mehdi Jalalmaab and Nasser L. Azad**

**Abstract**

In this study, a stochastic power management strategy for in-wheel motor electric vehicles (IWM-EV) is proposed to reduce the energy consumption and increase the driving range, considering the unpredictable nature of the driving power demand. A stochastic dynamic programming (SDP) approach, policy iteration algorithm, is used to create an infinite horizon problem formulation to calculate optimal power distribution policies for the vehicle. The developed SDP strategy distributes the demanded power, $P_{dem}$, between the front and rear IWMs by considering states of the vehicle, including the vehicle speed and the front and the rear wheels' slip ratios. In addition, a skid avoidance rule is added to the power management strategy to maintain the wheels' slip ratios within the desired values. Undesirable slip ratios cause poor brake and traction control performances and therefore should be avoided. The resulting strategy consists of a time-invariant, rule-based controller which is fast enough for real-time implementations, and additionally, it is not expensive to be launched since the future power demand is approximated without a need to vehicle communication systems or telemetric capability. A high-fidelity model of an IWM-EV is developed in the Autonomie/Simulink environment for evaluating the proposed strategy. The simulation results show that the proposed SDP strategy is more efficient in comparison to some benchmark strategies, such as an equal power distribution (ED) and generalized rule-based dynamic programming (GRDP). The simulation results of different driving scenarios for the considered IWM-EV shows the proposed power management strategy leads to 3% energy consumption reduction in average, at no additional cost. If the resulting energy savings is considered for the total annual trips for the vehicle, and also, the total number of electric vehicles in the country, the proposed power management strategy has a significant impact.

## Introduction

Battery electric vehicle (BEV) is a promising transportation technology with zero greenhouse gas emissions and no fossil fuel consumption that all the driving torque is delivered by electric motor(s). Electric motors are much more efficient than internal combustion engines and deliver torques faster and more accurately to the wheels. Recent advances in lithium-ion battery technologies and power electronics place BEVs in the center of many attentions.[1] In-wheel motor electric vehicle

(IWM-EV) is one type of BEVs which, despite of one central electric motor in conventional BEVs, is powered by two or four electric motors attached inside the wheels that work independently from each other in two operating modes; namely, driving and regenerative braking. Consequently, this architecture improves the controllability and power efficiency of the electric vehicle.

Performance, stability and maneuverability of IWM-EVs have been investigated in the literature.[2-8] The IWM-EV architecture increases the vehicle's maximum available yaw moment up to 40% and enhances the turning performance.[2,3] Additionally, the IWM technology improves the performance of motion control strategies e.g. dynamic traction force distribution and direct yaw moment control.[4–6] Interestingly, in-wheel motors can be employed for better estimations of many vehicle states and parameters, including the wheel slip ratios, slip angles and also body's slip angle due to an accurate output torque provided by this technology.[8] Jalali et al.[7] developed an advanced torque vectoring controller to generate the required corrective yaw moment through the torque intervention of the individual IWMs to stabilize the vehicle during both normal and emergency driving maneuvers. In addition to these advantages, the compact drivetrains inside the wheels free up space and allow the designers to optimize the vehicle layout and present new vehicle concepts. Figure 1 shows an IWM-EV's architecture and Figure 2 depicts IWM subsystems in more detail.

One prominent obstacle to wide applications of IWM-EVs and other BEVs is their low driving range and a slow charging process. Therefore, power management strategies that maximize power consumption efficiency receive high level of attention for BEVs. In this work, we propose a high-level power management strategy for IWM-EVs to determine proper split of the driver's demanded power between the vehicle motors/generators to maintain their operating points at high-efficiency regions, reducing the battery energy consumption while taking into account the vehicle system constraints. Traditional BEVs have no choice to manage the energy flow because they have only one traction motor which should generate all the demanded power for the wheels, but IWM-EVs have multiple power sources, hence there is degrees of freedom for the power management controller to minimize the electrical energy consumption.

Despite the hybrid electric vehicles (HEV) and plug-in hybrid electric vehicles (PHEV) power management problem which have been widely studied,[9–13] the IWM-EV power management research is relatively an open area. A good review of the current IWM-EV power management strategies is presented by Liu et al.[14] Power management based on IWM efficiencies have been investigated in several published works.[15–21] Qian et al.[22] proposed one time-step horizon simple search optimization algorithm for the front and rear IWMs torque distribution and achieved up to 27.4% reduction in the power consumption in comparison with a traditional EV. IWM-EV energy optimization using the Dynamic Programming (DP) algorithm based on terrain profile preview was investigated for constant speed driving situations[23] and a specific traffic model.[24] The results showed that the terrain profile changes the optimal torque distribution between the front and rear IWMs, and also, around 29% energy consumption improvement can be achieved. A nonlinear model predictive controller (NMPC) for the regenerative braking control of IWM-EVs is presented by Huang & Wang.[25] The proposed controller improves the energy recovery of the regenerative braking by determining the front and rear braking torques independently. Another

Department of Systems Design Engineering, University of Waterloo, Canada

**Corresponding author:**

Mehdi Jalalmaab, Department of Systems Design Engineering, University of Waterloo, 200 University Ave. W, Waterloo, N2L 3G1, Ontario, Canada.

Email: mjalalma@uwaterloo.ca

powerful technique for the IWM-EV power management problem is control allocation which was basically developed to distribute the desired total control effort among a redundant set of actuators for an over-actuated system.[17–21] Optimal control approaches using the Pontryagin's minimum principle are other candidates for online power management optimization.[26]

For deterministic driving scenarios and power demand profiles, DP is an excellent off-line technique to find the optimal power or torque distribution; but in the real-world problems, the future driving power demand is uncertain and dependent on the road traffic, terrain profile, weather, and so on. Moreover, the DP algorithm is computationally expensive and cannot be used for real-time applications. In this study, we will employ a Stochastic Dynamic Programming (SDP) algorithm to include the stochastic nature of power management for our real-time control problem. The proposed SPD strategy provides several advantages compared with other power management methods. First, this control scheme is cost-effective; in contrast to intelligent vehicle systems which employ extra sensors and communication systems to predict the vehicle speed profile and future load, SDP just processes the past driving information to predict the future driving condition without requiring any new sensors. Secondly, the SDP rules can be implemented as a time-invariant feedback algorithm or look-up table for real-time driving applications easily. Third, SDP can handle nonlinear cost functions and constraints seamlessly and we can take advantage of this potential to integrate safety-related and stability controllers, such as skid avoidance and antilock brake systems, to the power management strategy.

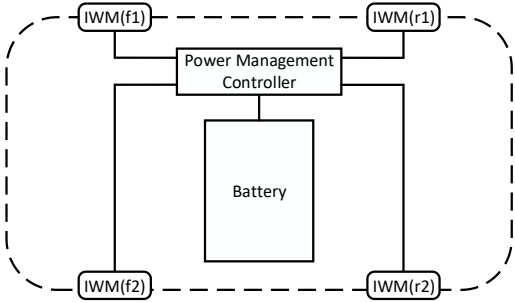

**Figure 1.** A schematic of IWM-EV architecture.

The SDP optimal control technique can handle constrained nonlinear optimization problems under uncertainties.[27] This technique has been used for an optimal energy management in HEVs [28,29] and PHEVs [30,31] to find the optimal power-split between the ICE and electric motors. For these platforms, it was shown that the power management control strategy using shortest path stochastic dynamic programming (SP-SDP) yields 2-3% better performance than the current factory-made controller for real-world driving situations.[29] Due to differences in the power management control problem as well as the vehicle system dynamics and constraints for IWM-EVs, we will investigate potential of an SDP-based power management strategy for this type of electric vehicle propulsions in this work.

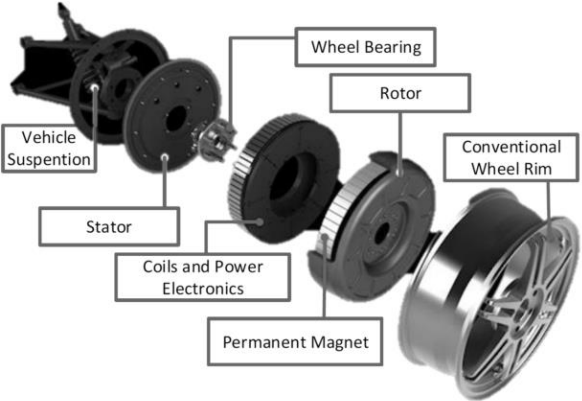

**Figure 2.** A typical in-wheel motor's components (courtesy by Protean).[38]

In the next sections, first, a high-fidelity model of an IWM-EV along with a control-oriented model to design the power management controller is introduced. Then, the stochastic power management strategy will be developed and the simulation results to evaluate the proposed algorithm using the high-fidelity model in the Autonomie/Simulink environment will be given. Finally, the main characteristics of the proposed strategy based on the simulation results are discussed.

## Modeling

In this section, two different IWM-EV models are developed for the design and the evaluation of proposed power management controller: a high-fidelity model and a control-oriented model. The high-fidelity model is a relatively accurate and detailed representation of IWM-EVs, including full models of the vehicle subsystems, which makes it suitable for a sensitivity analysis and the control system evaluation. On the other hand, the control-oriented model is a simpler and much faster model which is developed by ignoring low-level dynamic subsystems, simplifying complex equations, and eliminating unrelated vehicle component models to characterize the most prominent states of the vehicle system from the power management controller's point of view. The SDP control parameter optimization is an iterative convergence process, and therefore, a complex, computationally expensive high-fidelity model is not suitable for it. Thus, the control-oriented model is used for the control system parameter optimization and high-fidelity model is employed for the sensitivity analysis and control evaluation purposes. The vehicle specifications for the vehicle modeling are also listed in Table 1, obtained from Chen et. al. [23]

### High-fiedelity model development

In this study, the Autonomie software is used to develop the high-fidelity model. Autonomie is a modeling and simulation package developed by the Argonne National Lab, to develop automotive system models and evaluate vehicle control systems in a simulated environment. It includes many pre-built vehicle systems models, thus it allows users to simulate their desired models simply and efficiently.[32,33] Using Autonomie version-1210, a BEV model was modified to represent the desired IWM-EV. The developed high-fidelity vehicle model consists of a battery, IWMs, wheels and a chassis. The battery model contains 72 battery cells. Inside the battery cell block, the battery cell voltage, $V_{cell}$ is calculated by utilizing some look-up tables for SoC-$V_{cell}$, SoC-Ohmic and SoC-Polarization resistances. The IWM model determines the both maximum continuous and instantaneous torques based on the rotational speed and a heat index to saturate the commanded torques before applying them to the wheels.

**Table 1**. Baseline vehicle parameters.[23]

| Symbol | Parameters | Values |
| --- | --- | --- |
| $M$ | Vehicle mass | 800 kg |
| $E_{bat,max}$ | Battery capacity | 200 Ah |
| $V_{bat}$ | Nominal voltage of each battery cell | 3.3 V |
| $N_{bat}$ | Number of battery cells in series | 72 |
| $P_{motor}$ | In-wheel motor maximum power | 7.5 kw |
| $A$ | Vehicle front area | 1.66 m$^2$ |
| $\rho$ | Air density | 1.2 $\frac{kg}{m^3}$ |
| $L$ | Vehicle wheelbase | 1.84 m |
| $L_f$ | Front wheelbase | 0.92 m |
| $L_r$ | Rear wheelbase | 0.92 m |
| $h$ | Height of center of gravity | 0.6 m |
| $R_{eff}$ | Tire effective radius | 0.33 m |
| $T_{brake}^{regen,max}$ | Maximum regenerative braking torque | 80 Nm |
| $\mu_{max}$ | Maximum friction coefficient | 0.8 |

*Control-oriented model development*

The control-oriented model consists of simple models of the chassis, IWM, tires and battery, as follows:

*Chassis Model:* Based on the vehicle forces shown in Figure 3, the longitudinal acceleration is given by:

$$a_x = \frac{1}{m}(F_d - F_a - F_{rr} - F_g), \quad (1)$$

where, $F_d$ is the summation of all wheel driving forces, $F_a$ is the aerodynamic force, $F_{rr}$ is the total rolling resistance force, and $F_g$ is the road slope force which is equal to zero for a flat road driving:

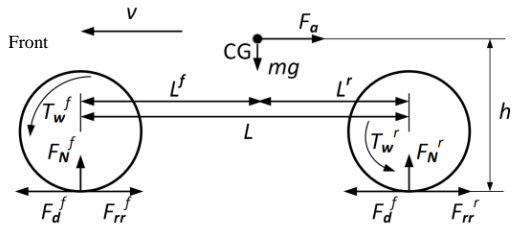

**Figure 3.** A schematic of the vehicle forces and moments.

$$F_d = \sum_{i=\{f1,f2,r1,r2\}} F_d^i, \quad (2)$$

$$F_a = \frac{1}{2}\rho A C_a v^2, \quad (3)$$
$$F_{rr} = mg \cos\theta\, C_{rr}, \quad (4)$$
$$F_g = mg \sin\theta, \quad (5)$$

where the superscript $i = \{f1, f2, r1, r2\}$ represents the wheel position: front $(f)$ or rear $(r)$ and left (1) or right (2). In the above equations, $\rho$ is the air density, $A$ is the effective area, and $C_a$ is the aerodynamic coefficient, $C_{rr}$ is the rolling resistance coefficient, $\theta$ is the slope of the road. The aerodynamic force, $F_a$, is a drag force caused by the movement of the vehicle in the presence of the air. The wheel driving forces are the main external forces that move the vehicle at the desired direction. They are a function of the normal force at the wheel and the friction coefficient:

$$F_d^i = F_N^i \mu^i, \quad (6)$$

whereas $\mu^i$ is the friction coefficient of wheel $i$ and $F_N^i$ is the wheel $i$ normal force. The normal forces can be determined by writing the moment equations around the contact points for the front and rear wheels. Referring to Figure 3, the moment equations can be written and rearranged, as follows:

$$F_N^{f1} = \frac{mg}{2L}(L^r \cos\theta - h \sin\theta) - \frac{1}{2L}(ma_xh + F_ah), \quad (7)$$
$$F_N^{r1} = \frac{mg}{2L}(L^f \cos\theta + h \sin\theta) + \frac{1}{2L}(ma_xh + F_ah). \quad (8)$$

*Tire Model:* The well-known, Magic tire formula[34], is used to determine the tire characteristics, for instance the friction coefficient. Based on the Magic formula, the wheel's friction coefficient is a nonlinear function of the slip ratio:

$$\mu_i = \mu_{max} D \sin(C \tan^{-1}(B\lambda_i - E(B\lambda_i - \tan^{-1}(B\lambda_i)))), \quad (9)$$

where $\mu_{max}$ is the maximum tire-road friction coefficient, $\lambda_i$ is the slip ratio of wheel $i$ and B, C,

D, and E are tire parameters assumed to be known.[25] The slip ratio is defined separately for driving and braking cases by:

$$\lambda_i = \frac{r\omega_i - v}{r\omega_i} \quad driving,$$
$$\lambda_i = \frac{r\omega_i - v}{v} \quad braking, \quad (10)$$

where $v$ is the vehicle speed and $\omega_i$ is the wheel rotational speed. In Figure 4, the friction coefficients versus slip ratio for the both high $\mu_{max}$ and low $\mu_{max}$ are shown.

In this study, only the longitudinal driving is considered, and therefore, the left and right motors will produce the same torques: $T_w^{f1} = T_w^{f2}$, $T_w^{r1} = T_w^{r2}$ and the angular velocities are the same for the right and left wheels all the time: $\omega_w^{f1} = \omega_w^{r2}$ and $\omega_w^{r1} = \omega_w^{r2}$. By taking into consideration all the moments of forces shown in Figure 3 about the center of front and rear wheels, the wheel rotation equations will be:

$$\dot{\omega}_w^i = \frac{1}{I_w^i}(T_w^i - F_d^i r_e^i), \quad (11)$$

whereas $I_w^i$ is the moment of inertia for the wheel $i$ and $r_e$ is the effective radius of the tire.

***In-Wheel Motor (IWM) Model:*** The demanded driving power, $P_{dem}$, is outlined by the driver through the throttle and brake pedals. A positive $P_{dem}$ is interpreted as the driving, and a negative $P_{dem}$ is considered as the brake command. For every motor/generator, two power flows are defined: $P_{MW}$, the power from (to) the motor (generator) to (from) the wheels and $P_{MB}$, the power from (to) the motor (generator) to (from) battery.

In the baseline IWM-EV, the motor/generators are directly mounted on the wheels, and thus, there are no transmission shafts, gear boxes and differentials, and consequently, there is no power loss between the motor/generators and wheels. Therefore, $P_{MW}$ is equal to the power exerted at the wheel $P_w$, i.e. $P_{MW} = P_w$. To satisfy the driver's commands, the sum of all $P_w$ should be equal to $P_{dem}$:

$$\sum_{i=\{f1,f2,r1,r2\}} P_w^i = P_{dem}. \quad (12)$$

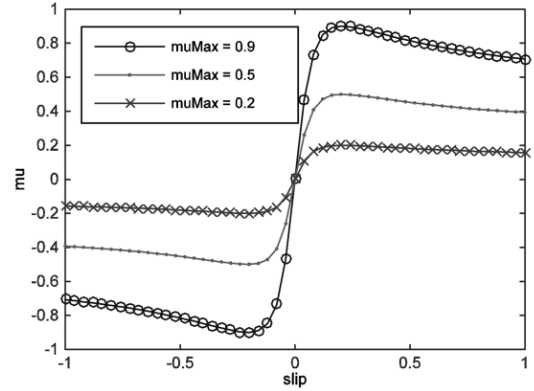

**Figure 4.** $\mu$-$\lambda$ relation by the Magic formula model.

Eq. (12) is called the drivability criteria. The power management controller specifies values of the motor powers $P_w^i$ by considering the drivability criteria as a constraint. Moreover, $P_{MB}$ equals to $P_{batt}$ which the battery is charging/discharging power:

$$P_{MB} = P_{batt} = I_{batt}V_{batt}, \quad (13)$$

where $V_{batt}$ and $I_{batt}$ are the battery discharging/charging voltage and current, respectively. Because of the frictional, thermal and other motor power losses, the motor-battery power $P_{MB}$ is not the same as the motor-wheel powers $P_{MW}$. The motor/generator efficiencies are defined as:

$$\eta_{trac} = \frac{P_{MW}}{P_{MB}} = \frac{P_w}{P_{batt}} \quad traction, \quad (14)$$

$$\eta_{regen} = \frac{P_{MB}}{P_{MW}} = \frac{P_{batt}}{P_w} \quad regen\ braking. \quad (15)$$

The motors' efficiency maps are illustrated in Figure 5. It can be concluded from these plots that efficiency maps are not strictly concave or convex. Also, when $\omega_w^i$ is low, the motor efficiency drops rapidly. The torques generated through the powers transferred to the wheels are calculated by:

$$T^i = \frac{P_w^i}{\omega^i} \qquad (16)$$

Similar to the torque and angular speed, in the straight-line driving, the powers at left and right wheels are the same:

$$P_w^{f1} = P_w^{f2} = P_w^f/2, \qquad (17)$$

$$P_w^{r1} = P_w^{r2} = P_w^r/2. \qquad (18)$$

Every motor has a maximum limit for the output power ($P_{w,max}^i$) and cannot produce a power more than that value.

$$P_w^i \leq P_{w,max}^i. \qquad (19)$$

For braking situations, the regenerative braking torque, $T_{brake}^{regen}$, has a certain limit, and thus, it should be constrained in the vehicle model. Consequently, at high braking torque demands, the friction brakes should be also involved to support the regenerative braking system to produce the required braking torques:

$$T_{brake} = T_{brake}^{fric} + T_{brake}^{regen}. \qquad (20)$$

In conclusion, some of the braking power cannot be harvested to charge the battery, and therefore, it is lost by the friction brakes.

*Battery:* The battery SoC equation is derived by utilizing a static equivalent circuit model (Figure 6). The battery current is given by:

$$I = \frac{V_{oc} + \sqrt{V_{oc}^2 - 4R_{batt}P_{batt}}}{2R_{batt}}, \qquad (21)$$

whereas $R_{batt}$ is the battery resistance and $P_{batt}$ is the battery power. By convention, a positive $P_{batt}$ indicates a discharging case and a negative $P_{batt}$ indicates a charging case. The rate of change of SoC for the battery is:

$$\frac{d}{dt}(SoC) = \frac{I}{E_{max}}, \qquad (22)$$

whereas $E_{max}$ is the total amount of charge that can be stored in the battery pack. The battery power is derived from the summation of all IWMs' consumed power:

$$P_{batt} = \sum_{i=\{f1,f2,r1,r2\}} P_w^i \cdot (\eta^i)^{-sgn(P_w^i)}, \qquad (23)$$

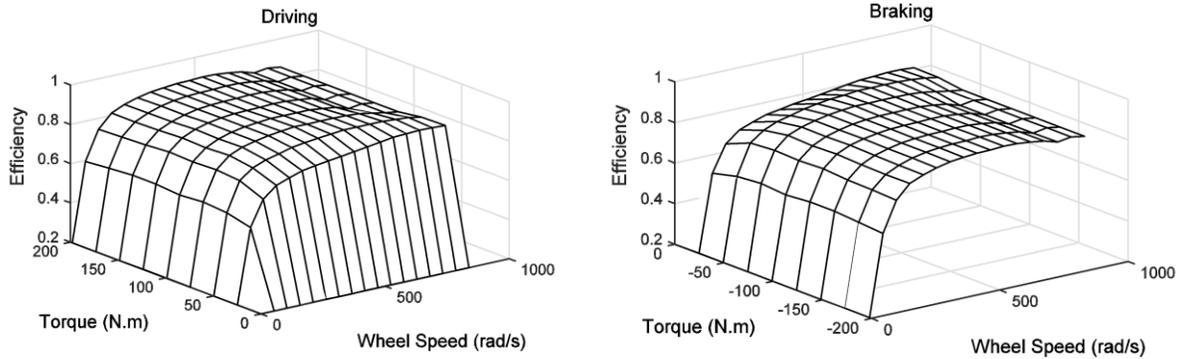

**Figure 5.** The motor efficiency maps for driving and braking (taken from the Autonomie software).

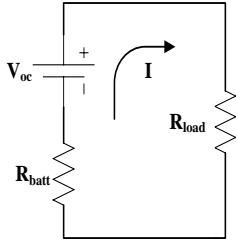

**Figure 6.** A battery equivalent circuit model.

where $\eta^i$ is the efficiency of the $i$-th IWM. The battery has a maximum limit for the power output during draining as well as the power input during recharging:

$$-P_{batt,max} \leq P_{batt} \leq P_{batt,max}. \qquad (24)$$

*Parameter estimation and model validation*

The high-fidelity model of the vehicle is used to estimate the key parameters of the control-oriented model, and also, validate the accuracy of the resulting model such that the simulation results of the both control-oriented model and Autonomie-based high-fidelity model will be close for the same drive cycle. In this way, the battery's resistance or $R_{batt}$ is estimated based on the simulation results for FTP-75 and HWFET drive cycles. To do so, the root mean square error given below is minimized:

$$RMSE = \sqrt{\sum_{i=0}^{N}(SoC_{HF}(i) - SoC_{CO}(i,R_{batt}))^2}, \qquad (25)$$

where $N$ is the number of samples and $SoC_{HF}$ and $SoC_{CO}$ are the instantaneous SoC of the high-fidelity model and the control-oriented model, respectively. By the minimization of RMSE for $R_{batt}$, we obtain:

$$R_{batt,est} = 0.063 \, \Omega. \qquad (26)$$

Figure 7 compares the states of control-oriented model and high-fidelity Autonomie model for FTP 75 drive cycle. The demanded power is distributed equally between the front and rear in-wheel motors. Figure 7 shows good agreement between the two models, and therefore, it can be concluded that the control-oriented model despite its simplicity, represents the key characteristics of the vehicle system properly.

*A sensitivity analysis of IWM-EV power management problem*

Before the controller design, we will conduct a sensitivity analysis to realize how the vehicle's states affect the power management problem. Minimizing the charge depletion or maximizing the driving range is the main goal of the power management controller; therefore, the power consumption variation in terms of the battery state of charge changes or $\Delta SoC$ is examined for a few scenarios in this section. The input of the power management controller is the front wheels power, $P_w^f$. More precisely, the effects of some vehicle's states on the optimal $P_w^f$ for one time step is studied. In Figures 8, 9 and 10 the plots of $\Delta SoC$ versus $P_w^f$ for different speed, slip ratios and $P_{dem}$ are shown.

In Figure 8, the variations of $P_w^f$ when the vehicle speed changes are depicted. In the low speed scenario, $v = 5 \, m/s$, $\Delta SoC$ minimization by a search method has a symmetrical solution for $P_w^f$. The same is true for the high speed scenario, $v = 20 \, m/s$. However, the optimal solution for $v = 10 \, m/s$ is unique and it is half of $P_{dem}$. Figure 9 shows the plot of $\Delta SoC$ versus $P_w^f$ for different front wheel slip ratios, $\lambda_f$, at the constant speed of 10 m/s and $P_{dem} = 10 \, kW$. We see that, for low slip ratios, $\lambda_f = 0.01$ and $\lambda_f = 0.2$, the optimal solution is almost symmetrical, but for a higher slip ratio $\lambda_f = 0.5$, the optimum $P_w^f$ is less than the half, $P_w^f = 3000 \, W$. Figure 10 demonstrates the optimal solution for different $P_{dem}$. We can see that, for all of the cases of Figure 10, the optimum solution is symmetrical.

As a result, the optimal power distribution for the IWMs is dependent to the vehicle's demanded power, speed and tires' slip rations.

## Stochastic power management controller design

In this section, a stochastic model of the driver's power demand is developed, and then, the SDP algorithm is used to solve the real-time power management problem of IWM-EV. The SDP algorithm or the Markov decision process is a computational technique for solving stochastic, state-dependent optimization problems to find the optimum sequences of control actions.

It is assumed that $P_{\text{dem}}$ is a Markovian variable, that is, the probability distribution of the next step $P_{dem}$ is just a function of the current time-step states (transition probability). By assuming the Markov property, the power management problem can be solved by a SDP algorithm to determine the optimal power distribution between the front and rear wheels of IWM-EV.

Since the vehicle dynamics equation and the power management cost function are assumed to be time-invariant, and additionally, no final time cost term or constraint is defined for the cost function, the problem is formulated as an infinite horizon problem[35]. The infinite horizon formulation generates a set of the time-invariant control policies, which can be easily used for real-time control applications. The approximate policy iteration algorithm is also used to solve the infinite horizon SDP problem[36]. The proposed stochastic power management strategy employs the expected total discounted reward. For the infinite horizon

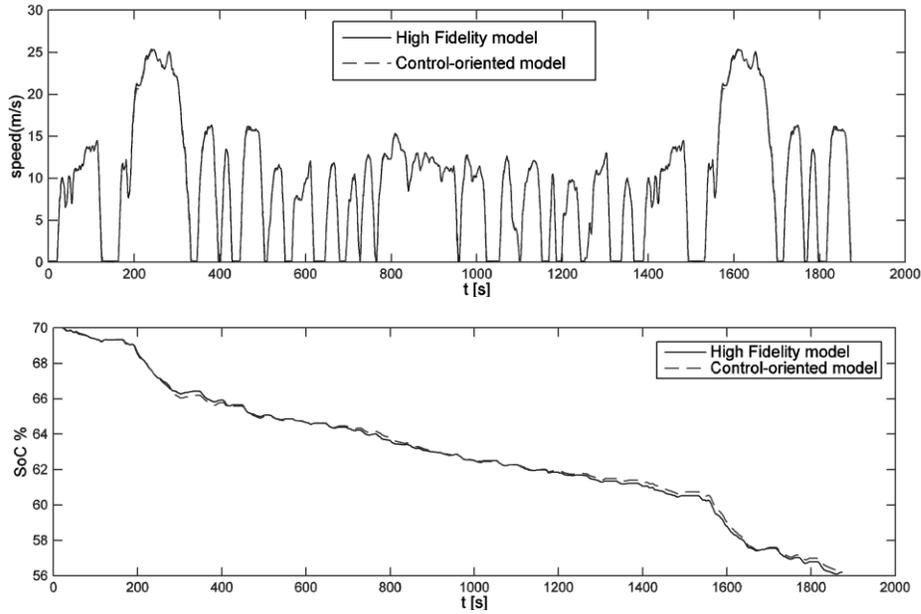

**Figure 7.** Comparison of speed and SoC for the Autonomie and control-oriented mode for the FTP 75 drive cycle.

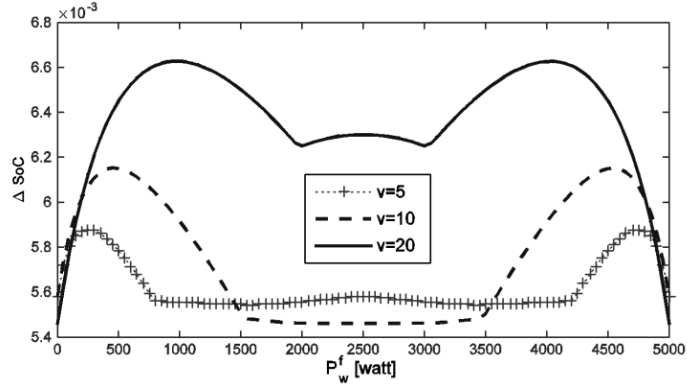

**Figure 8.** The effect of vehicle speed on $\Delta SoC$.

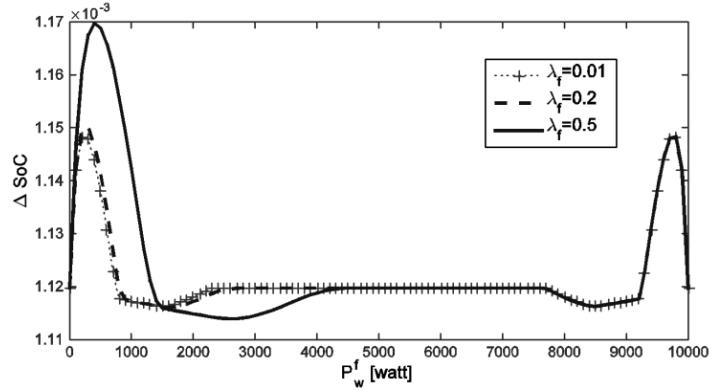

**Figure 9.** The effect of $\lambda_f$ on $\Delta SoC$.

Markovian SDP problem, the expected total discounted cost has the following form:

$$J_\pi(x_0) = \lim_{N \to \infty} \mathop{E}_{w_k} \left\{ \sum_{k=0}^{N-1} \gamma^k \xi(x_k, \pi(x_k)) \right\}, \quad (27)$$

where $J_\pi(x_0)$ is the expected cost, $\pi(x)$ is the control policy, $\xi$ is the one-time step cost, and $\gamma$ is the discount factor, $0 < \gamma < 1$. For the IWM-EV power management problem, the cost is a sum of the battery SoC changes and a penalty for $P_{dem}$ deviation, as given below:

$$\xi = \Delta SoC + \alpha. M, \quad (28)$$

whereas $M$ is the power error and $\alpha$ is a weighting factor. $M$ is measured by a squared error between the demanded power and the actual power produced by the wheels:

$$M = (P_{dem} - P_f - P_r)^2, \quad (29)$$

The implemented SDP algorithm (policy iteration) has two main calculation steps: the policy evaluation step and the policy improvement step. At each iteration, first, the cost function $J_\pi(x)$ is calculated for the given policy (policy evaluation):

$$J_\pi^{s+1}(x^i) = \xi\left(x^i, \pi(x^i)\right) + \mathop{E}_{\{P_{dem(i+1)}\}} \{\gamma J_\pi^s(\acute{x})\}, \quad (30)$$

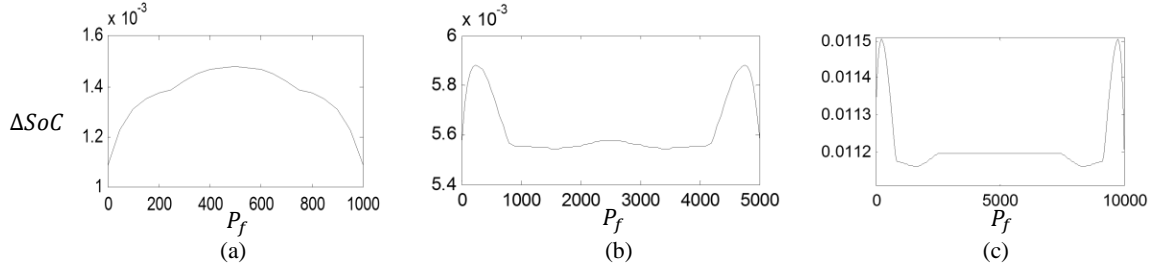

**Figure 10.** The effects of power distribution on $\Delta SoC$ for: (a) $P_{dem} = 1\ kw$, (b) $P_{dem} = 5\ kw$ and (c) $P_{dem} = 10\ kw$.

where $\acute{x}$ shows the propagated states at the end of time step, $\acute{x} = f\left(x^i, \pi(x^i)\right)$. Afterwards, by using the calculated $J_\pi$, the policy is updated by minimizing the following equation for each element of the state space, $i$:

$$\pi(x^i) = \underset{u \in U(x^i)}{argmin}\left[\xi(x^i, u) + \underset{P_{dem(i+1)}}{E}\{\gamma J_\pi(\acute{x})\}\right], \quad (31)$$

The optimized policy of this step is returned to the step one to update the cost function again. This process is repeated until $\pi(x)$ converges and does not change considerably over the next iterations. In Figure 11, a flowchart of the policy iteration algorithm is indicated.

Figure 12 depicts the norm of changes in $\pi(x)$ over subsequent iterations and shows how this norm converges to zero at the end of the algorithm, which means $\pi(x)$ does not change anymore.

A problem definition with a proper number of the states is critical for the SDP algorithm since this algorithm develops some rule policies as a function of the entire state sets, and consequently, a large number of the states will make the problem solving process slow and computationally expensive, which is considered as the curse of dimensionality. Thus, the least number of the states should be chosen to capture the dominant dynamics of the system. First, it is worth mentioning that, despite the PHEV/HEV power management problems, $SoC$ is not a considerable state for the IWM-EV power management problem and it does not affect the optimal power distribution calculation. This is mainly due to the fact that the battery is the only source of the power, and thus, regardless of the battery charge level, it should be used to provide the demanded power.

In the sensitivity analysis section, the effects of $\{P_{dem}, v, \lambda_f, \lambda_r\}$ on the optimal power management solution were explained. As shown in Equation (22), the rate of $\Delta SoC$ is a function of $P_{batt}$ and $P_{batt}$ changes with the motor efficiencies, $\eta^i_{trac}$ and $\eta^i_{regen}$, and also, the demanded power, $P_{dem}$, as shown in Equations (14) and (15). Considering Equations (9) and (10), it can be stated that $\Delta SoC$ is a function of $P_{dem}, v$ and tires slip ratios. Accordingly, the states of SDP problem for this study are:

$$x = [P_{dem}, v, \lambda_f, \lambda_r]^T, \quad (32)$$

Since the slip ratio is only a function of the vehicle speed and the wheel rotational speed, the rotational speeds are substituted by the slip ratio in the problem state sets. By doing this, we make sure that the state space will cover all of the slip ratio values which are critical for the effective braking/accelerations in harsh driving situations.

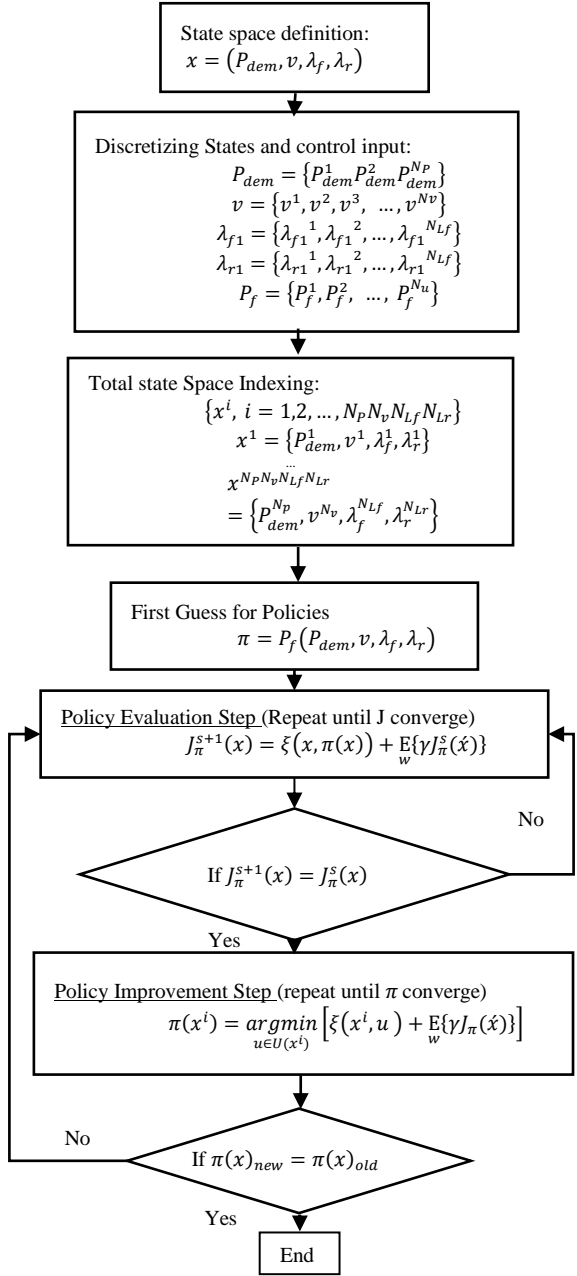

**Figure 11.** Policy iteration algorithm flowchart.[36]

To solve the SDP problem, it is necessary to discretize the state domain of the problem. The demanded power $P_{dem}$ is assumed to take a finite number of values, where $N_P$ is the number of $P_{dem}$ discretization:

$$P_{dem} \in \{P_{dem}^1, P_{dem}^2, \dots, P_{dem}^{N_P}\}. \tag{33}$$

Other states are discretized as:

$$\lambda_{f1} \in \{\lambda_{f1}^{\ 1}, \lambda_{f1}^{\ 2}, \dots, \lambda_{f1}^{\ N_{Lf}}\} \tag{34}$$

$$\lambda_{r1} \in \{\lambda_{r1}^{\ 1}, \lambda_{r1}^{\ 2}, \dots, \lambda_{r1}^{\ N_{Lf}}\}, \tag{35}$$

$$v \in \{v^1, v^2, \dots, v^{N_v}\}, \tag{36}$$

Then, the total state space is:

$$\{x^i, i = 1, 2, \dots, N_P N_v N_{Lf} N_{Lr}\}. \tag{37}$$

The front in-wheel motor power $P_f$ is considered as the control input. This variable is also discretized as $\{P_f^1, P_f^2, \dots, P_f^{N_u}\}$.

An appropriate discretization resolution of the power demand as well as the time span of the decision period should be simultaneously chosen such that a proper transition probability matrix (TPM) of $P_{dem}$ can be created. By a trial-and-error, we found the discretization sets for the problem states as given below:

$P = -12 : 1 : 19\ kW$,

$v = [0\ 5\ 10\ 25]\ m/s$,

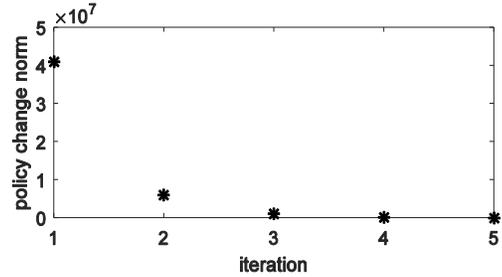

**Figure 12.** The policy iteration algorithm convergence behavior.

$\lambda_f = \lambda_r = [0, \pm 0.001, \pm 0.1, \pm 0.21, \pm 0.35, \pm 1]$,

$dt_{SDP} = 0.1 \ sec$.

Since $P_{dem}$ is assumed to be a Markovian state, we will have:

$$Pr\{P_{dem,k+1} = P_{dem}^j | P_{dem,k} = P_{dem}^i\} = p_{i,j}, \quad (38)$$
$$\forall i, j = 1, 2, \ldots, N_P$$

where, the summation of all next step probabilities for a specific state value should be:

$$\sum_{j=1}^{N_P} p_{i,j} = 1, \quad (39)$$

where $p_{i,j}$ is a one-step transition probability when the system demanded power is $P_{dem}^j$ at time $k+1$, and $P_{dem}^i$ at time $k$.

### Constraints

The first constraint is the driving torque limit of the electric motors:

$$T_i \leq T_{max\ motor} \quad (for\ driving), \quad (40)$$

The next constraint is related to the regenerative braking. During the vehicle braking, if the demanded brake torque exceeds the limit, it is divided between the regenerative and frictional braking, as follows:

$$T_i = T_{regen} + T_{friction}, \quad (41)$$

$$T_{regen} \leq T_{max\ regen} \quad (for\ breaking), \quad (42)$$

The last constraint is about the bounds on the battery delivered power:

$$P_{batt,min} \leq P_{batt} \leq P_{batt,max}, \quad (43)$$

### Skid Avoidance System Integration with Power Management Strategy

To prevent safety problems while using the power management strategy, skid avoidance constraints are applied as a part of the power management control strategy of the IWM-EV. In the conventional skid avoidance and anti-lock brake (ABS) systems, a high frequency bang–bang controller in the form of an on–off controller is used to keep the slip ratios in the desired interval to maximize the braking performance. Since the electric motor's response time is much faster than typical ABS of ICE vehicles,[37] the task of ABS can be easily handled by electric motor's power controller.

In this study, a set of skid avoidance constraints are proposed to be considered in the development of the stochastic dynamic programing policies. The proposed constraints for the slips are presented in Table 2. The main logic behind proposing these constrains is that, in each IWM, while performing hard brakes, if the wheel's slip ratio violates the desired interval, the power management controller should set the brake command to zero to relax the braking system and avoid the violation. The critical slip ratios are equal to $\pm 0.2$, and thus the desired slip ratio interval is: $\lambda = [-0.2, +0.2]$. Therefore, constraints of Table 2 are applied. When $\lambda_f < -0.2$, the front IWM commanded power, $P_f$, should be zero; and the same is for the rear IWM. Consequently, the only situation in which the power distribution problem is over-actuated and has two degrees of freedom is the case that the both wheels' slip ratio are inside the desired interval.

**Table 2.** Skid Avoidance Constraints.

|  | $\lambda_f < -0.2$ | $\lambda_f \geq -0.2$ |
|---|---|---|
| $\lambda_r < -0.2$ | $P_f = 0$<br>$P_r = 0$ | $P_f = P_{dem}$<br>$P_r = 0$ |
| $\lambda_r \geq -0.2$ | $P_f = 0$<br>$P_r = P_{dem}$ | Search for the best distribution |

## Simulation Results

The SDP strategy is optimized using the control-oriented vehicle model. The output of the policy iteration algorithm is an SDP policy set, that is, a set of the control rules specifying the front and rear wheels power distribution on the basis of the current states of the vehicle. After determining the control policies, we need to evaluate the optimality of the strategy by comparing it with some other techniques.

The utilized SDP technique requires adequate observations to build the TPM matrix, which is necessary for the policy optimizations. Three driving cycles, namely FTP, HWFET and NYCC, have been used as the observation driving scenarios to create TPM. Also, to evaluate the resulting stochastic power management strategy, we simulated its performance for the drive cycles used for observations (FTP, HWFET, and NYCC), and additionally a new drive cycle UDDS, to investigate the robustness of the devised power management strategy for new driving demands.

Furthermore, to compare the performance of the proposed SDP controller, we performed some simulations for two other power management strategies. The first is front-rear equal distribution (ED) of all power demands, and the second is the Generalized Rule-based DP (GRDP) strategy, as the benchmarks for the comparison. ED is a fixed ratio distribution policy which distributes the power demands equally and independent of the system's states between the front and the rear wheels in all instances. Our preliminary investigation shows that an equal distribution is the most optimal fixed ratio for the IWM-EV power management problem.

Since the DP solution is drive cycle dependent, the rules from one drive cycle may cause undesirable results for another drive cycle, thus, the GRDP method [35] is proposed as the second benchmark. To determine the GRDP strategy, first, the DP-based problem is solved for a specific drive cycle, that is, FTP 75, and then, we obtain a simple linear rule for the power distribution as a function of $P_{dem}$. The DP solution for the FTP driving cycle and a single rule extracted from $P_f$ versus $P_{dem}$ profile is shown in Figure 13.

Tables 3, 4 and 5 show the comparison of $\Delta SoC$ for different drive cycles and road conditions based on the high-fidelity model simulation. A lower $\Delta SoC$ is more preferred since it means that less power is consumed, and thus, the power management strategy is more efficient and the vehicle's driving range is longer. The last columns express the improvement percentage for the SDP strategy in comparison with the ED strategy. Although the DP strategy results in the global optimal solutions, but a generalization of the DP results to develop the GRDP policies deteriorate the optimality of these rules for different drive cycles, and they are, generally, worse than the ED and SDP rules.

By comparing these tables, it is concluded that, for all road conditions, SDP works better than ED and GRDP. The maximum improvement is achieved in the NYCC drive cycle with a slippery road condition ($\mu_{max} = 0.2$), resulting in a 22% reduction in the power consumption. From these tables, we can see that, the SDP strategy's power efficiency increases by a decrease of $\mu_{max}$, which means that the SDP strategy successfully considers the road friction and the tire slip ratio in the decision making process.

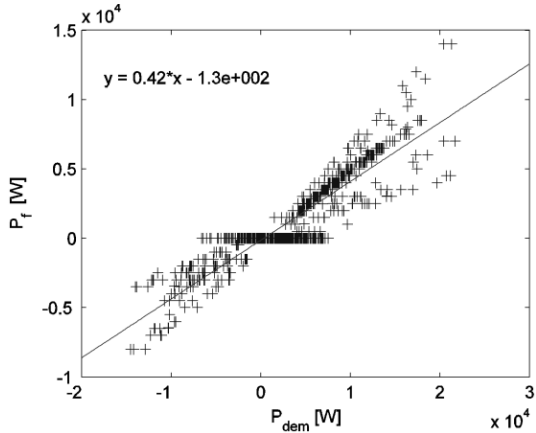

**Figure 13.** $P_f$ versus $P_{dem}$ for the DP solution for the FTP 75 drive cycle.

To evaluate the performance of the SDP strategy for new drive cycles which were not used in the SDP policy development, the results for the UDDS drive cycle in Table 3, Table 4 and Table 5 show that SDP improves the power efficiency for this drive cycle. Therefore, we can conclude that the SDP approach provides a robust power management strategy for various drive cycles in real-time applications.

Some parts of $P_f$ and $P_r$ variation curves versus the time for the UDDS drive cycle over a snowy road ($\mu_{max} = 0.2$) are shown in Figure 14. The whole UDDS drive cycle takes 1370 sec, and thus only, a part of the variations for $t = [345\ 430]$ is shown in this figure. As it was stated before, ED assigns a half of the demanded power to the front IWMs and another half to the rear IWMs. The Figure 14 shows that the SDP strategy assigns more than a half of the power to the front wheels most of the time. In particular, in the time interval of [370,380], the difference between SDP and ED is considerable and $P_r$ for the SDP strategy in this interval is almost zero. Due to the IWMs efficiency maps, this decision prevents having power loss from the rear IWMs, and also, makes the front IWMs to generate power at a more efficient working point.

The slip ratios of the front and rear tires for the same road ($\mu_{max} = 0.2$) are depicted in Figure 15. Without the skid avoidance constraints, the slip ratio of the front wheel converges to -1 in many braking instances, which means that the wheel is completely locked. The proposed power management with skid avoidance constraints optimizes the power consumption and keeps the wheels' slip ratios in the desired intervals simultaneously. From Figure 15, it can be seen that this strategy has maintained the braking slip rations in the acceptable interval successfully. Based on the results given in Table 6, adding the skid avoidance constraints to the stochastic power management strategy does not change the power depletion of the vehicles.

**Table 3.** $\Delta SoC$ for simulations with $\mu_{max} = 0.9$.

| Drive Cycle | ED | GRDP | SDP | %$\Delta SoC$ (SDP vs. ED) |
|---|---|---|---|---|
| FTP | 13.8 | 13.79 | 13.71 | 0.65 |
| HWFET | 14.88 | 14.86 | 14.68 | 1.34 |
| NYCC | 1.585 | 1.58 | 1.58 | 0.32 |
| UDDS | 8.78 | 8.77 | 8.72 | 0.68 |

**Table 4.** $\Delta SoC$ for simulation with $\mu_{max} = 0.5$.

| Drive Cycle | Fixed Ratio | DP Rules | SDP | %$\Delta SoC$ (SDP vs. ED) |
|---|---|---|---|---|
| FTP | 13.89 | 13.87 | 13.6 | 2.09 |
| HWFET | 14.92 | 14.90 | 14.73 | 1.27 |
| NYCC | 1.6 | 1.59 | 1.53 | 4.37 |
| UDDS | 8.84 | 8.82 | 8.67 | 1.92 |

**Table 5.** $\Delta SoC$ for simulation with $\mu_{max} = 0.2$.

| Drive Cycle | Fixed Ratio | DP Rules | SDP | %$\Delta SoC$ (SDP vs. ED) |
|---|---|---|---|---|
| FTP | 14.3 | 14.29 | 14.19 | 0.77 |
| HWFET | 15.04 | 15.01 | 14.93 | 0.73 |
| NYCC | 3.41 | 3.4 | 2.64 | 22.58 |
| UDDS | 9.09 | 9.09 | 9.07 | 0.22 |

## Conclusions

The contributions of this investigation are listed below:
- Development of a SDP-based power management strategy for IWM-EVs to consider the stochastic nature of driving commands.
- Incorporation of wheel slip ratios into the IWM-EV power management strategy to distribute power demands optimally in different road conditions.
- Adding skid avoidance constraints to the power management controller prevents the wheels from

locking in severe braking while maintains the optimality of the power management strategy.

The evaluation of the proposed power management strategy for different driving cycles and road friction coefficients shows that it can decrease the battery charge depletion, while providing the demanded power for the IWMs. In average, the SDP strategy saved 3% of SoC in comparison with an equal energy distribution strategy. As a result, the vehicle employing the SDP-based power management strategy can follow the considered drive cycles properly, and eventually, it has a higher final SoC as compared to the other strategies, such as ED and GRDP.

Although the design and evaluation procedures were carried out for a certain baseline vehicle, the algorithms presented in this work are general and can be implemented for other IWM-EVs. The power distribution between the left and right wheels has been considered equal. Considering curvy roads requires an unequal distribution of the power between the left and right IWMs which adds another level of complexity to the problem in hand. This can be considered as future work for this study.

**Table 6.** $\Delta SoC$ for simulations with $\mu_{max} = 0.9$.

| Drive Cycle | SDP | SDP with skid avoidance |
|---|---|---|
| FTP | 13.71 | 0.65 |
| HWFET | 14.68 | 1.34 |
| NYCC | 1.58 | 0.32 |
| UDDS | 8.72 | 0.68 |

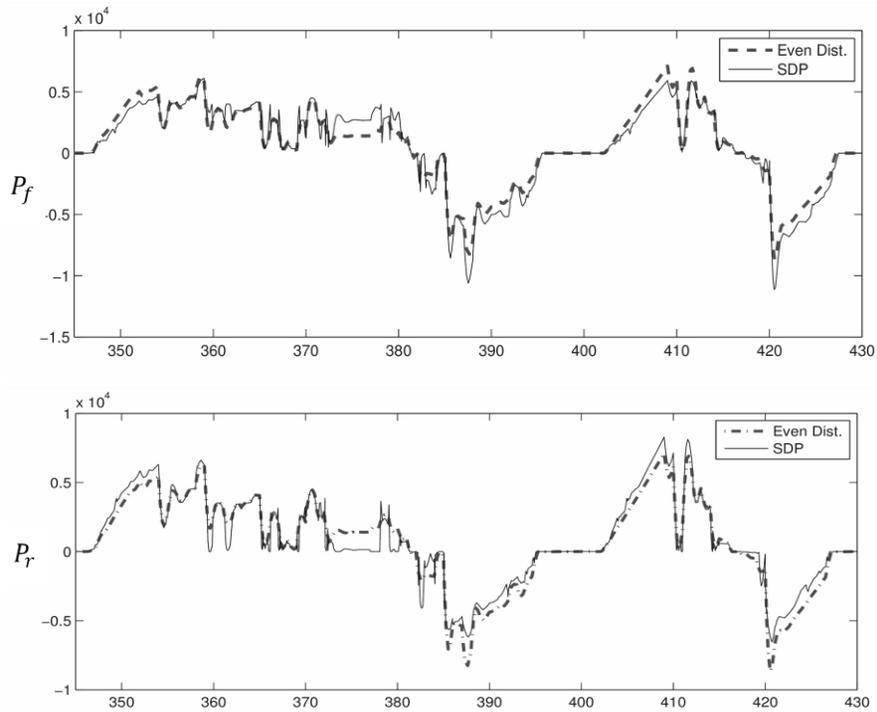

**Figure 14.** Variations of $P_f$ and $P_r$ for SDP and ED for the UDDS drive cycle with low friction coefficient.

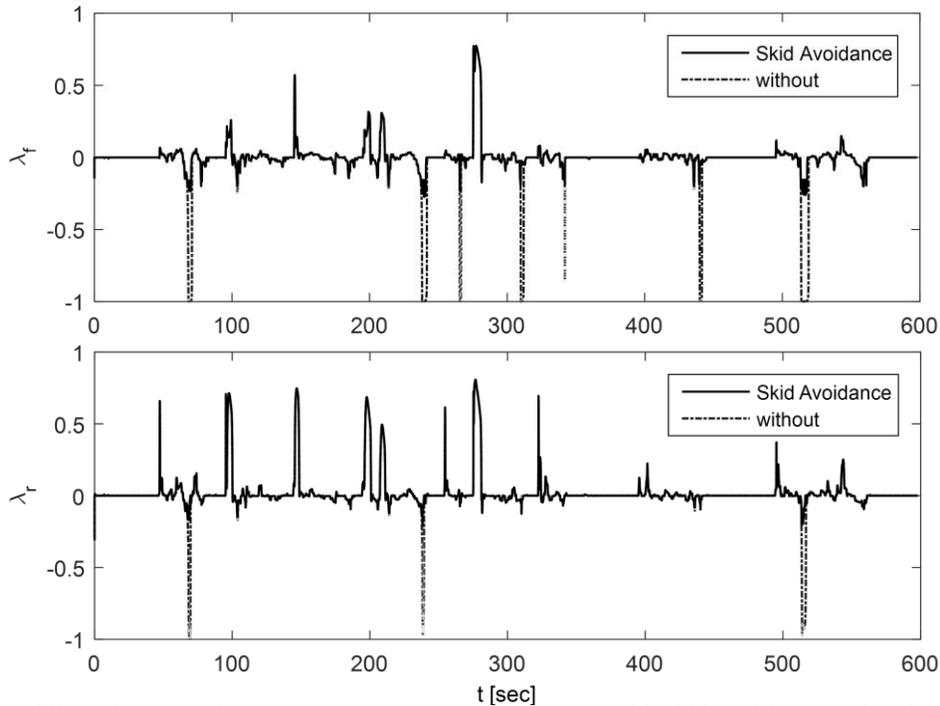

**Figure 15.** Slip ratio comparison for power management systems with skid avoidance and without skid avoidance constraint for (a) front and (b) rear wheels.

**Acknowledgements**

The authors would like to thank Maryyeh Chehresaz for her help in developing the vehicle high-fidelity model in the Autonomie environment.

**References**

1. Trigg T, Telleen P. *Global EV Outlook*. 2013.
2. Jacobsen B. Potential of electric wheel motors as new chassis actuators for vehicle manoeuvring. *Proc Inst Mech Eng Part D J Automob Eng* 2002; 216: 631–640.
3. Kim W, Yi K. Coordinated Control of Tractive and Braking Forces Using High Slip for Improved Turning Performance of an Electric Vehicle Equipped with In-Wheel Motors. In: *IEEE Vehicular Technology Conference*. 2012.
4. Sakai SI, Sado H, Hori Y. Motion control in an electric vehicle with four independently driven in-wheel motors. *IEEE/ASME Trans Mechatronics* 1999; 4: 9–16.
5. Chang C. Torque distribution control for electric vehicle based on traction force observer. In: *IEEE International Conference on Computer Science and Automation Engineering*. 2011, pp. 371–375.
6. He P, Hori Y, Kamachi M, et al. Future motion control to be realized by in-wheel motored electric vehicle. In: *Annual Conference of IEEE Industrial Electronics Society*. 2005, pp. 2632–2637.
7. Jalali K, Uchida T, Lambert S, et al. Development of

an Advanced Torque Vectoring Control System for an Electric Vehicle with In-Wheel Motors using Soft Computing Techniques. *SAE Int J Altern Powertrains* 2013; 2: 261–278.
8. Geng C, Uchida T, Hori Y. Body slip angle estimation and control for electric vehicle with in-wheel motors. In: *33rd Annual Conference of the IEEE Industrial Electronics Society*. 2007, pp. 351–355.
9. Vajedi M, Chehresaz M, Azad NL. Intelligent Power Management of Plug-in Hybrid Electric Vehicles, Part I: Real-time Optimum SOC Trajectory Builder. *Int J Electr Hybrid Veh* 2013; 6: 46–67.
10. Vajedi M, Chehresaz M, N. L. Azad 2013. Intelligent Power Management of Plug-in Hybrid Electric Vehicles, Part II: Real-time Route- based Power Management. *Int J Electr Hybrid Veh* 2013; 6: 68–86.
11. Taghavipour A, Azad NL, McPhee J. Design and evaluation of a predictive powertrain control system for a plug-in hybrid electric vehicle to improve the fuel economy and the emissions. *Proc Inst Mech Eng Part D J Automob Eng* 2015; 229: 624–640.
12. Taghavipour A, Azad NL, McPhee J. An Optimal Power Management Strategy for Power Split Plug-in Hybrid Electric Vehicles. *Int J Veh Des* 2012; 60: 286–304.
13. Razavian RS, Taghavipour A, Azad NL, et al. Design and Evaluation of a Real-Time Optimal Control System for Series Hybrid Electric Vehicles. *Int J Electr Hybrid Veh* 2012; 4: 260–288.
14. Liu H, Chen X, Wang X. Overview and Prospects on Distributed Drive Electric Vehicles and Its Energy Saving Strategy. *Prz Elektrotechniczny* 2012; 122–125.
15. Xu W, Zheng H, Liu Z. The Regenerative Braking Control Strategy of Four-Wheel-Drive Electric Vehicle Based on Power Generation Efficiency of Motors. In: *SAE World Congress & Exhibition*. 2013.
16. Gu J, Ouyang M, Lu D, et al. Energy efficiency optimization of electric vehicle driven by in-wheel motors. *Int J Automot Technol* 2013; 14: 763–772.
17. Chen Y, Wang J. Adaptive Energy-Efficient Control Allocation for Planar Motion Control of Over-Actuated Electric Ground Vehicles. *IEEE Trans Control Syst Technol* 2014; 22: 1362–1373.
18. Laine L, Fredriksson J. Coordination of Vehicle Motion and Energy Management Control Systems for Wheel Motor Driven Vehicles. In: *IEEE Intelligent Vehicles Symposium*. 2007, pp. 773–780.
19. Chen Y, Wang J. A global optimization algorithm for energy-efficient control allocation of over-actuated systems. In: *American Control Conference*. 2011, pp. 5300–5305.
20. Chen Y, Wang J. Energy-efficient control allocation with applications on planar motion control of electric ground vehicles. In: *American Control Conference*. 2011, pp. 2719–2724.
21. Chen Y, Wang J. Fast and Global Optimal Energy-Efficient Control Allocation With Applications to Over-Actuated Electric Ground Vehicles. *IEEE Trans Control Syst Technol* 2011; 1202–1211.
22. Qian H, Xu G, Yan J, et al. Energy management for four-wheel independent driving vehicle. In: *IEEE International Conference on Intelligent Robots and Systems*. 2010, pp. 5532–5537.
23. Chen Y, Chen Z, Wang J. Operational Energy Optimization for Pure Electric Ground Vehicles Based on Terrain Profile Preview. In: *ASME Dynamic Systems and Control Conference*. 2011, pp. 271–278.
24. Chen Y, Li X, Wiet C, et al. Energy management and driving strategy for in-wheel motor electric ground vehicles with terrain profile preview. *IEEE Trans Ind Informatics* 2014; 10: 1938–1947.
25. Huang X, Wang J. Model predictive regenerative braking control for lightweight electric vehicles with in-wheel motors. *Proc Inst Mech Eng Part D J Automob Eng* 2012; 226: 1220–1232.
26. Azad NL. Online optimization of automotive engine coldstart hydrocarbon emissions control at idle conditions. *Proc Inst Mech Eng Part I J Syst Control Eng* 2015; 229: 781–796.
27. Bertsekas DP. *Dynamic Programming and Optimal Control: Volume II*. Athena Scientific, 2007.
28. Lin C-C, Peng H, Grizzle JW. A stochastic control strategy for hybrid electric vehicles. In: *Proceedings of the American Control Conference*. 2004, pp. 4710–4715.
29. Opila DF, Xiaoyong Wang, McGee R, et al. Performance comparison of hybrid vehicle energy management controllers on real-world drive cycle data. In: *American Control Conference*. 2009, pp. 4618–4625.
30. Moura SJ, Fathy HK, Callaway DS, et al. A Stochastic Optimal Control Approach for Power Management in Plug-In Hybrid Electric Vehicles. *IEEE Trans Control Syst Technol* 2010; 19: 545–555.
31. Zhang H, Qin Y, Li X, et al. Driver-Oriented Optimization of Power Management in Plug-In Hybrid Electric Vehicles. In: *EIC Climate Change*

*Technology Conference*. 2013, pp. 1–12.
32. Vijayagopal R, Michaels L, Rousseau AP, et al. Automated Model Based Design Process to Evaluate Advanced Component Technologies. In: *SAE World Congress & Exhibition*. 2010.
33. Kim N, Rousseau A. Comparison between rule-based and instantaneous optimization for a single-mode, power-split HEV. In: *SAE World Congress & Exhibition*. 2011.
34. Pacejka HB, Bakker E. The Magic Formula Tyre Model. *Veh Syst Dyn* 1992; 21: 1–18.
35. Grizzle JWJW, Lin C, Peng H, et al. Power management strategy for a parallel hybrid electric truck. *IEEE Trans Control Syst Technol* 2003; 11: 839–849.
36. Puterman ML. *Markov Decision Processes: Discrete Stochastic Dynamic Programming*. John Wiley & Sons, 2009.
37. Fujimoto H. Regenerative Brake and Slip Angle Control of Electric Vehicle with In-wheel Motor and Active Front Steering. In: *Proc. 1st Internatinal Electric Vehicle Technology …*. 2011.
38. PROTEAN ELECTRIC. Protean Drive PD18 Subsystems.http://www.proteanelectric.com/en/subsystems/.